\documentclass{scrartcl}

\usepackage[T1]{fontenc}
\usepackage[utf8]{inputenc}
\usepackage{amsmath, amsthm, amssymb}
\usepackage{dsfont}
\usepackage{graphicx}
\usepackage[format=plain]{caption}
\usepackage{capt-of}
\usepackage{subcaption}
\usepackage{bm}
\usepackage{tikz}
\usetikzlibrary{arrows,shapes,trees,calc,through}

\newcommand{\norm}[1]{\left\| #1 \right\|}  
\renewcommand{\d}{\,\mathrm{d}} 
\newcommand{\e}{\mathrm{e}} 
\renewcommand{\i}{\mathrm{i}} 
\newcommand{\N}{\mathbb{N}}  
\newcommand{\Z}{\mathbb{Z}}

\newcommand{\R}{\mathbb{R}}
\newcommand{\C}{\mathbb{C}}

\renewcommand{\phi}{\varphi}

\renewcommand{\Re}{\operatorname{Re}}
\renewcommand{\Im}{\operatorname{Im}}

\numberwithin{equation}{section}

\newtheorem{thm}{Theorem}[section]

\newtheorem{prop}[thm]{Proposition}

\theoremstyle{definition} \newtheorem{ex}[thm]{Example}
\theoremstyle{definition} \newtheorem{df}[thm]{Definition}
\theoremstyle{definition}

\title{Infinite-time admissibility under compact perturbations}

\author{Jochen Schmid$^{1,2}$\\  
\small $^1$ Institut f\"ur Mathematik, Universit\"at W\"urzburg, 97074 W\"urzburg, Germany\\
\small $^2$ Fraunhofer Institute for Industrial Mathematics (ITWM), 67663 Kaiserslautern, Germany\\ 
\small jochen.schmid@itwm.fraunhofer.de}  

\date{}

\begin{document}

\maketitle

\begin{abstract}
\small{ \noindent 
We investigate the behavior of infinite-time admissibility under compact perturbations. 
We show, by means of two completely different examples, that infinite-time admissibility is not preserved under compact perturbations $Q$ of the underlying semigroup generator $A$, even if $A$ and $A+Q$ both generate strongly stable semigroups.  
}
\end{abstract}


\section{Introduction}

In this note, we investigate the behavior of infinite-time admissibility under compact perturbations of the underlying semigroup generator. So, we consider semigroup generators $A: D(A) \subset X \to X$ (with $X$ a Hilbert space) and possibly unbounded control operators $B$ (defined on another Hilbert space $U$) and we ask how the property of infinite-time admissibility of $B$ behaves under compact perturbations of the generator $A$. 
%
Infinite-time admissibility of $B$ for $A$ means that for every control input $u \in L^2([0,\infty),U)$ the mild solution of the initial value problem 
\begin{align} \label{eq:ivp-x_0=0}
x' = Ax + Bu(t) \qquad \text{and} \qquad x(0) = 0
\end{align}  
is a bounded 
function from $[0,\infty)$ with values in $X$. (A priori, the mild solution has values only in the extrapolation space $X_{-1}$ of $A$ and, a fortiori, need not be bounded in the norm of $X$, of course.) 
\smallskip

It is well-known that (finite-time) admissibility is preserved under very general perturbations $Q$ of the generator $A$, in particular, under bounded perturbations. It is also clear that infinite-time admissibility, by contrast, is not preserved under bounded perturbations. 
Just think of a generator $A$ of an exponentially stable semigroup and a bounded perturbation $Q$ (for example, a sufficiently large multiple of the identity) such that $A+Q$ has 
spectral points in the right half-plane.
\smallskip

In this note, we will show by way of two completely different kinds of examples that infinite-time is also not preserved under compact perturbations $Q$ which are such that both $A$ and $A+Q$ generate stronly stable (but not exponentially stable) semigroups. So, in other words, we show that there exist semigroup generators $A$ and $A+Q$ with $Q$ being compact and a control operator $B$ such that
\begin{itemize}
\item the semigroups $\e^{A\cdot}$ and $\e^{(A+Q)\cdot}$ are strongly stable but not exponentially stable
\item $B$ is infinite-time admissible for $A$ but not infinite-time admissible for $A+Q$.
\end{itemize} 
In our first -- more elementary -- example, we will use an old and well-known result from the 1970s, namely a stabilization result for collocated linear systems. In that example, the compact perturbation $Q$ will be of rank $1$ and the control operator $B$ will be bounded. In particular, none of the technicalities 
coming along with unbounded control operators will bother us there. 
In our second -- less elementary -- example, we will use a more advanced result from the 1990s, namely a characterization of infinite-time admissibility for diagonal semigroup generators. In that example, the control operator $B$ will be unbounded and the compact perturbation $Q$ will be of rank $\infty$.
\smallskip

In the entire note, we will use the following notation. 
\begin{align*}
\R^+_0 := [0,\infty), \qquad \C^+ := \{z \in \C: \Re z > 0\}, \qquad
\C^- := \{z \in \C: \Re z < 0\}.
\end{align*}
As usual, $L(X,Y)$ denotes the Banach space of bounded linear operators between two Banach spaces $X$ and $Y$ and $\norm{\cdot}_{X,Y}$ stands for the operator norm on $L(X,Y)$. Also, $\norm{u}_2$ denotes the norm of a square-integrable function $u \in L^2(\R^+_0,U)$ with values in the Banach space $U$. And finally, for a semigroup generator $A$ and bounded operators $B, C$ between appropriate spaces, the symbol $\mathfrak{S}(A,B,C)$ will stand for the state-linear system~\cite{CuZw} 
\begin{align*}
x' = Ax + Bu(t) \quad \text{with} \quad y(t) = Cx(t).
\end{align*}

\section{Some basic facts about admissibility and infinite-time admissibility}

In this section, we briefly recall the definition of and some basic facts about admissibility and infinite-time admissibility. 
If $A: D(A) \subset X \to X$ is a semigroup generator on the Hilbert space $X$ and $X_{-1}$ is the corresponding extrapolation space, then an operator $B \in L(U,X_{-1})$ (with $U$ another Hilbert space) is called \emph{control operator for $A$}. Also, $B$ is called a \emph{bounded control operator} iff $B \in L(U,X)$ and an \emph{unbounded control operator} iff $B \in L(U,X_{-1}) \setminus L(U,X)$. See ~\cite{TuWe} (Section~2.10) or~\cite{EnNa} (Section~II.5) for basic facts about extrapolation spaces.  

\begin{df}
Suppose $A: D(A) \subset X \to X$ is a semigroup generator on $X$ and $B \in L(U,X_{-1})$, where $X, U$ are both Hilbert spaces. Then $B$ is called \emph{admissible for $A$} iff for every $u \in L^2(\R^+_0,U)$
\begin{align} \label{eq:adm-def}
(0,\infty) \ni t \mapsto \Phi_t(u) := \int_0^t \e^{A_{-1}s} B u(s) \d s
\end{align}
is  a function with values in $X$, where $A_{-1}$ is the generator of the continuous extension of the semigroup $\e^{A\cdot}$ to $X_{-1}$. 
\end{df}

Clearly, for a given semigroup generator $A$ every bounded control operator $B \in L(U,X)$ is admissible (because $\e^{A_{-1}s}|_X = \e^{As}$ for $s \in \R^+_0$).  It should also be noted that if $B \in L(U,X_{-1})$ is admissible for $A$, then for every $t \in (0,\infty)$ the linear operator $L^2(\R^+_0,U) \ni u \mapsto \Phi_t(u) \in X$ defined in~\eqref{eq:adm-def} is closed and thus continuous  by the closed graph theorem. Consequently, $B \in L(U,X_{-1})$ is admissible for $A$ if and only if
\begin{align}
\Phi_t \in L(L^2(\R^+_0,U),X) \qquad (t \in (0,\infty)).
\end{align}

\begin{df}
Suppose $A: D(A) \subset X \to X$ is a semigroup generator on $X$ and $B \in L(U,X_{-1})$, where $X, U$ are both Hilbert spaces. Then $B$ is called \emph{infinite-time admissible for $A$} iff for every $u \in L^2(\R^+_0,U)$
\begin{align} \label{eq:infin-time-adm-def}
(0,\infty) \ni t \mapsto \Phi_t(u) := \int_0^t \e^{A_{-1}s} B u(s) \d s
\end{align}
is a function with values in $X$ that is bounded (in the norm of $X$), where $A_{-1}$ is the generator of the continuous extension of the semigroup $\e^{A\cdot}$ to $X_{-1}$. 
\end{df}

Clearly, if $B \in L(U,X_{-1})$ is infinite-time admissible for a given semigroup generator $A$, then it is also admissible for $A$. It should also be noted that, by the uniform boundedness principle, $B \in L(U,X_{-1})$ is infinite-time admissible for $A$ if and only if
\begin{align} \label{eq:equiv def infin-time adm}
\Phi_t \in L(L^2(\R^+_0,U),X) \qquad (t \in (0,\infty)) 
\qquad \text{and} \qquad
\sup_{t\in (0,\infty)} \norm{\Phi_t}_{L^2(\R^+_0,U),X} < \infty.
\end{align}
Some authors~\cite{Oo00}, \cite{CuWe06}, \cite{WeCu06} use the term input-stability for the system $\mathfrak{S}(A,B)$ instead of infinite-time admissibility.
\smallskip

It is well-known that admissibility is preserved under bounded perturbations.

\begin{prop}
Suppose $A: D(A) \subset X \to X$ is a semigroup generator on $X$ and $B \in L(U,X_{-1})$, where $X, U$ are both Hilbert spaces. Also, let $Q \in L(X)$. Then $B$ is admissible for $A$ if and only if $B$ is admissible for $A+Q$. 
\end{prop}

In fact, the conclusion of this proposition remains true for much more general perturbations $Q$, namely for perturbations of the (feedback) form $Q = B_0C_0$, where $B_0 \in L(U_0,X_{-1})$ is an admissible control operator for $A$ and $C_0 \in L(X,U_0)$ with $U_0$ an arbitrary Hilbert space. See Corollary~.5.5.1 from~\cite{TuWe}, for instance.

\begin{prop}
Suppose $A: D(A) \subset X \to X$ is the generator of an exponentially stable semigroup on $X$ and $B \in L(U,X_{-1})$ is admissible for $A$. Then $B$ is even infinite-time admissible for $A$. 
\end{prop}

See Proposition~4.4.5 in~\cite{TuWe}, for instance, and notice that for bounded control operators $B$ the  above proposition is trivial. In view of that proposition, it is clear that infinite-time admissibility -- unlike admissibility -- is not preserved under bounded perturbations. Choose, for example, a bounded generator $A$ of an exponentially stable semigroup and let $Q := -A \in L(X)$ and $B := I \in L(X,X)$ (identity operator on $X$). 

\section{An example using a stabilization result for collocated linear systems}

\subsection{Stabilization of collocated linear systems}

We will use the following well-known stabilization result for collocated systems, that is, systems of the form $\mathfrak{S}(A,B,B^*)$  with a bounded control operator $B$. It essentially goes back to~\cite{Be78} (Corollary~3.1) and, in the form below, can be found in~\cite{Oo00} (Lemma~2.2.6), for instance. (Actually, for the more general version with  the countability assumption on $\sigma(A_0) \cap \i \R$ we have to refer to~\cite{WeCu06}, but this more general version will not be used in the sequel.)

\begin{thm} \label{thm:stabiliz colloc lin syst}
Suppose $A_0$ is a contraction semigroup generator on a Hilbert space $X$ with compact resolvent (or, more generally, with $\sigma(A_0) \cap \i \R$ being countable). Suppose further that $B \in L(U,X)$ with another Hilbert space $U$ and that $\mathfrak{S}(A_0,B,B^*)$ is approximately controllable in infinite time or approximately observable in infinite time. Then 
\begin{itemize}
\item[(i)] $B$ is infinite-time admissible for $A_0-BB^*$, more precisely, 
\begin{align*}
\norm{ \int_0^t \e^{(A_0-BB^*)s} B u(s) \d s}_X^2 \le \frac{1}{2} \norm{u}_2^2
\qquad (u \in L^2(\R^+_0,U), t \in \R^+_0).
\end{align*}
\item[(ii)] $\e^{(A_0-BB^*)\cdot}$ is a strongly stable contraction semigroup on $X$. 
\end{itemize}
\end{thm}

A far-reaching generalization of this result to the case of unbounded control operators was obtained by Curtain and Weiss~\cite{WeCu06}. See Theorem~5.1 and 5.2  in conjunction with Proposition~1.5 from~\cite{WeCu06}. 
We also refer to~\cite{CuWe06} for a parallel result on exponential stabilization. 

\subsection{Infinite-time admissibility under compact perturbations}

\begin{ex}
Set $X := \ell^2(\N,\C)$ and let $A_0: D(A_0) \subset X \to X$ be defined by
\begin{align*}
A_0 x := (\lambda_{0 k} x_k)_{k \in \N} \qquad (x \in D(A_0)),
\end{align*}
where $D(A_0) := \{(x_k) \in X: (\lambda_{0 k} x_k) \in X \}$ and $\lambda_{0 k} := -\alpha_k + \i \beta_k$ with
\begin{align*}
\Re \lambda_{0 k} = -\alpha_k := -1/k \qquad (k \in \N) 
\qquad \text{and} \qquad 
\Im \lambda_{0 k} = \beta_k \longrightarrow \infty \qquad (k \to \infty).
\end{align*}
Set $U := \C$ and let $B: U \to \C^{\N}$ be defined by 
\begin{align*}
B u := (u b_k)_{k\in\N} \qquad (u \in U),
\end{align*}
where
\begin{align*}
b_k := 1/k^{3/8} \qquad (k\in I_1) \qquad &\text{and} \qquad b_k := 1/k \qquad (k\in I_2) \\
I_1 := \{ l^2: l\in \N \} \qquad &\text{and} \qquad I_2 := \N \setminus I_1.
\end{align*}
Clearly, $(b_k) \in X$ and therefore $B \in L(U,X)$. We now define 
\begin{align*}
A := A_0 - BB^* \qquad \text{and} \qquad A' := A_0
\end{align*}
and show, in various steps, that $A$ and $A'$ are generators of strongly but not exponentially stable contraction semigroups on $X$, that $A' = A+Q$ for a compact perturbation $Q$ of rank one, and that $B$ is infinite-time admissible for $A$ but not infinite-time admissible for $A'$. 
\smallskip

As a first step, we observe that $A' = A+Q$ with $Q := BB^*$ and that $Q$ has rank one (because the same is true for $B$), whence $Q$ is compact. 
\smallskip

As a second step, we observe from
\begin{align*}
\lambda_{0 k} \in \C^- \qquad (k \in \N) \qquad \text{and} \qquad \sup\{ \Re \lambda_{0 k}: k \in \N\} = 0
\end{align*}
that $A'$ is the generator of a strongly stable but not exponentially stable contraction semigroup on $X$. 
\smallskip

As a third step, we show that $B$ is not infinite-time admissible for $A'$. 
In view of~\eqref{eq:equiv def infin-time adm} we have to show that 
\begin{align} \label{eq:ex-1, 1}
\sup_{\norm{u}_2=1} \sup_{t\in (0,\infty)} \norm{ \int_0^t \e^{A_0 s} Bu(s) \d s }_X = \infty.
\end{align}
We first observe by Fatou's lemma that 
\begin{align} \label{eq:ex-1, 2}
\liminf_{t\to\infty} \norm{ \int_0^t \e^{A_0 s} Bu(s) \d s }_X^2 
&\ge \sum_{k\in\N} \bigg| \int_0^{\infty} u(s) \e^{\lambda_{0 k} s} \d s \bigg|^2 |b_k|^2 \notag \\
&\ge \bigg| \int_0^{\infty} u(s) \e^{\lambda_{0 n} s} \d s \bigg|^2 |b_n|^2
\end{align}
for every $u \in L^2(\R^+_0,U)$ and $n \in \N$. Setting $u_n(s) := n^{-1/2} \chi_{[0,n]}(s) \cdot \e^{-\i \beta_n s}$ for $s \in \R^+_0$ and $n \in \N$, we see that
\begin{align}
\norm{u_n}_2 &= 1 \label{eq:ex-1, 3}\\
\bigg| \int_0^{\infty} u_n(s) \e^{\lambda_{0 n} s} \d s \bigg|^2 
= \frac{1}{n} \bigg| \int_0^n \e^{-\alpha_n s} \d s \bigg|^2 
&= \frac{1}{\alpha_n^2 n} \Big( 1- \e^{-\alpha_n n} \Big)^2 = n (1-\e^{-1})^2 \label{eq:ex-1, 4}
\end{align}
for every $n \in \N$. Combining now~\eqref{eq:ex-1, 2}, \eqref{eq:ex-1, 3} and~\eqref{eq:ex-1, 4} we get
\begin{align*}
\sup_{\norm{u}_2=1} \sup_{t\in (0,\infty)} \norm{ \int_0^t \e^{A_0 s} Bu(s) \d s }_X^2 
&\ge \sup_{n\in\N} \bigg( \liminf_{t\to\infty} \norm{ \int_0^t \e^{A_0 s} Bu_n(s) \d s }_X^2 \bigg) \notag \\
&\ge  (1-\e^{-1})^2 \sup_{n\in\N} \big( n |b_n|^2 \big).
\end{align*}
Since $\sup_{n\in\N} \big( n |b_n|^2 \big) \ge \sup_{n\in I_1} \big( n |b_n|^2 \big)  = \infty$, 
the desired relation~\eqref{eq:ex-1, 1} follows.
\smallskip

As a fourth step, we show that $B$ is infinite-time admissible for $A$ and that $A$ is the generator of a strongly stable contraction semigroup on $X$. 
In order to do so, we apply the stablization theorem above (Theorem~\ref{thm:stabiliz colloc lin syst}). Since 
\begin{align*}
\Re \lambda_{0 k} \le 0 \qquad (k \in \N) \qquad \text{and} \qquad |\lambda_{0 k}| \longrightarrow \infty \qquad (k \to \infty),
\end{align*}
we see that $A_0$ is a contraction semigroup generator on $X$ with compact resolvent, and since the eigenvalues $\lambda_{0 k}$ of $A_0$ are pairwise distinct and $b_k \ne 0$ for every $k \in \N$, we see that the collocated linear system $\mathfrak{S}(A,B,B^*)$ is approximately controllable and approximately observable in infinite time (Theorem~4.2.3 of~\cite{CuZw}). So, by the stablization theorem above (Theorem~\ref{thm:stabiliz colloc lin syst}), $B$ is infinite-time admissible for $A_0-BB^* = A$ and $\e^{A\cdot}$ is a  strongly stable contraction semigroup on $X$. 
\smallskip

As a fifth and last step, we convince ourselves that the semigroup generated by $A$ is not exponentially stable. 
Assume the contrary. Then there exist $M \ge 1$ and $\omega < 0$ such that $\{z \in \C: \Re z > \omega \} \subset \rho(A)$ and
\begin{align*}
\norm{(A-z)^{-1}} \le \frac{M}{\Re z - \omega} \qquad (\Re z > \omega).
\end{align*} 
So, since $\Re \lambda_{0 n} \longrightarrow 0$ as $n \to \infty$, we conclude that 
\begin{align} \label{eq:ex-1, 5}
\limsup_{n\to\infty} \norm{(A-\lambda_{0 n})^{-1}} \le \limsup_{n\to\infty}  \frac{M}{\Re \lambda_{0 n} - \omega} = \frac{M}{|\omega|}.
\end{align}
We now observe that
\begin{align} \label{eq:ex-1, 6}
(A-\lambda_{0 n}) e_n = -BB^* e_n = -b_n \cdot b \longrightarrow 0 \qquad (n \to \infty).
\end{align}
Combining~\eqref{eq:ex-1, 5} and~\eqref{eq:ex-1, 6} we arrive at
\begin{align*}
1 = \limsup_{n\to\infty} \norm{e_n} = \limsup_{n\to\infty} \norm{ (A-\lambda_{0 n})^{-1} b_n \cdot b } \le  \frac{M}{|\omega|} \limsup_{n\to\infty} \norm{ b_n \cdot b } = 0.
\end{align*}
Contradiction! $\blacktriangleleft$
\end{ex}

\section{An example using an admissibility result for diagonal linear systems}

\subsection{Characterization of infinite-time admissibility}

We will use the following well-known characterization of infinite-time admissibility for diagonal semigroup generators $A_0$. It essentially goes back to~\cite{Gr95} (Proposition~2.2) and can also be found in~\cite{TuWe} (Theorem~5.3.9 in conjunction with Remark~4.6.5), for instance.

\begin{thm} \label{thm:char infinite-time adm}
Suppose $X = \ell^2(I,\C)$ with a countable infinite index set $I$ and let $A_0: D(A_0) \subset X \to X$ be the diagonal operator given by
\begin{align*}
A_0 x := (\lambda_{0k}x_k)_{k\in I} \qquad (x \in D(A_0)),
\end{align*}
where $D(A_0) := \{(x_k) \in X: (\lambda_{0k}x_k) \in X\}$ and $\lambda_{0k} \in \C^-$ for every $k \in I$. Suppose further that $B \in L(U,X_{-1})$ with $U := \C$, that is, 
\begin{align*}
Bu = (u b_k)_{k\in I} \qquad (u \in U)
\end{align*}
for a uniqe sequence $(b_k) \in X_{-1} = \{ (c_k) \in \C^{I}: \sum_{k\in I} |c_k|^2/(1+|\lambda_k|^2) < \infty\}$. Then the following statements are equivalent:
\begin{itemize}
\item[(i)] $B$ is infinite-time admissible for $A_0$
\item[(ii)] there exists a constant $M \in \R_0^+$ such that 
\begin{align*}
\sum_{k\in I} \frac{|b_k|^2}{|z-\lambda_{0k}|^2} \le \frac{M}{\Re z} \qquad (z \in \C^+).
\end{align*}
\end{itemize}
\end{thm}

Clearly, in the situation of the above theorem the condition~(ii) is equivalent to the existence of  a constant $M \in \R_0^+$ such that 
\begin{align} \label{eq:resolv-estim-infin-time-adm}
\norm{(z-A)^{-1}B}_{U,X} \le \frac{M}{\sqrt{\Re z}} \qquad (z \in \C^+).
\end{align}
A far-reaching generalization of the above theorem to the case of general contraction semigroup generators $A_0$ on a separable Hilbert space $X$ was obtained by Jacob and Partington~\cite{JaPa01}. See Theorem~1.3 from~\cite{JaPa01}. It states that for a contraction semigroup generator $A_0$ on a separable Hilbert space $X$ a control operator $B \in L(U,X_{-1})$ with $U := \C$ is infinite-time admissible if and only if there is a constant $M \in \R_0^+$ such that the resolvent estimate~\eqref{eq:resolv-estim-infin-time-adm} is satisfied. 
We also refer to~\cite{JaPa04} and~\cite{TuWe} (Section~5.6) for an overview of many more admissibility results, for example, for infinite-dimensional input-value spaces~$U$.

\subsection{Infinite-time admissibility under compact perturbations}

\begin{ex}
Set $X := \ell^2(\Z,\C)$ and let $A: D(A) \subset X \to X$ and $A': D(A') \subset X \to X$ be defined by
\begin{align*}
Ax := (\lambda_k x_k)_{k\in \Z} \qquad (x \in D(A)) 
\qquad \text{and} \qquad
A'x := (\lambda_k' x_k)_{k\in \Z} \qquad (x \in D(A')),
\end{align*}
where $D(A) := \{(x_k) \in X: (\lambda_k x_k) \in X\}$ and $D(A') := \{(x_k) \in X: (\lambda_k' x_k) \in X\}$ with
\begin{align*}
\lambda_k := \begin{cases} -1/k^{1/2} + \i k, \quad k \in \N \\ -(|k|+1)^{1/2}, \quad k \in -\N_0 \end{cases}
\qquad \text{and} \qquad
\lambda_k' := \begin{cases} -\e^{-k} + \i k, \quad k \in \N \\ -(|k|+1)^{1/2}, \quad k \in -\N_0 \end{cases}.
\end{align*}
Set $U := \C$ and let $B: U \to \C^{\Z}$ be defined by
\begin{align*}
B u := (u b_k)_{k\in \Z} \qquad (u \in U),
\end{align*}
where
\begin{align*}
b_k := 1/k \qquad (k \in \N), \qquad b_0 := 0, \qquad b_k := 1/|k|^{1/2} \qquad (k \in -\N).
\end{align*}
Clearly, $\sum_{k\in\Z}  |b_k|^2/(1+|\lambda_k|^2) < \infty$ and $\sum_{k\in \Z} |b_k|^2 = \infty$ whence 
$(b_k) \in X_{-1} \setminus X$. And therefore
\begin{align*}
B \in L(U,X_{-1}) \setminus L(U,X).
\end{align*}
We now show, in various steps, that $A$ and $A'$ are generators of strongly but not exponentially stable contraction semigroups on $X$, that $A' = A+Q$ for a compact perturbation $Q$ of infinite rank, and that $B$ is infinite-time admissible for $A$ but not infinite-time admissible for $A'$. 
\smallskip

As a first step, we observe from
\begin{align*}
\lambda_k, \lambda_k' \in \C^- \qquad (k \in \Z) 
\qquad \text{and} \qquad
\sup \{\Re \lambda_k : k \in \Z\}, \sup \{\Re \lambda_k' : k \in \Z\} = 0
\end{align*}
that $A$ and $A'$ are generators of strongly stable but not exponentially stable contraction semigroups on $X$. 
\smallskip

As a second step, we observe that $A' = A+Q$ for a compact operator $Q$ of infinite rank.
Indeed, the operator $Q: X \to X$ defined by
\begin{align*}
Qx := ((\lambda_k'-\lambda_k)x_k)_{k\in\Z} \qquad (x \in X)
\end{align*}
is a bounded operator on $X$ because $(\lambda_k'-\lambda_k)_{k\in\Z}$ is a bounded sequence. Also, $Q$ is the limit in norm operator topology 
of the finite-rank operators $Q_N: X \to X$ defined by
\begin{align*}
Q_N x := (\dots, 0,0,(\lambda_1'-\lambda_1)x_1, \dots, (\lambda_N'-\lambda_N)x_N, 0, 0, \dots) 
\qquad (x \in X)
\end{align*} 
and therefore $Q$ is compact, as desired.
\smallskip

As a third step, we show that $B$ is infinite-time admissible for $A$. 
We have that 
\begin{align} \label{eq:ex-2, 1}
\sum_{k\in \Z} \frac{|b_k|^2}{|z-\lambda_k|^2} \le \sum_{k\in \Z} \frac{|b_k|^2}{(\Re z + |\Re \lambda_k|)^2}
\le \frac{1}{2 \Re z} \sum_{k\in \Z} \frac{|b_k|^2}{|\Re \lambda_k|} 
\end{align}
for every $z \in \C^+$ and that
\begin{align} \label{eq:ex-2, 2}
M:= \sum_{k\in \Z} \frac{|b_k|^2}{|\Re \lambda_k|}  < \infty.
\end{align}
So, by the admissibility theorem above (Theorem~\ref{thm:char infinite-time adm}), the claimed  infinite-time admissibility of $B$ for $A$ follows from~\eqref{eq:ex-2, 1} and~\eqref{eq:ex-2, 2}.
\smallskip

As a fourth and last step, we show that $B$ is not infinite-time admissible for $A'$. 
We have that 
\begin{align} \label{eq:ex-2, 3}
\sum_{k\in \Z} \frac{|b_k|^2}{|z-\lambda_k'|^2} 
\ge \frac{|b_n|^2}{|z-\lambda_n'|^2} = \frac{1}{(\Re z + \e^{-n})^2 + (\Im z - n)^2} \frac{1}{n^2}
\end{align}
for every $z \in \C^+$ and $n \in \N$. Choosing $z_n := \e^{-n} + \i n \in \C^+$ for $n \in \N$, we see from~\eqref{eq:ex-2, 3} that
\begin{align}
\sup_{z \in \C^+} \bigg( \Re z \sum_{k\in \Z} \frac{|b_k|^2}{|z-\lambda_k'|^2} \bigg) 
&\ge \sup_{n\in\N} \bigg( \frac{\Re z_n}{(\Re z_n + \e^{-n})^2 + (\Im z_n - n)^2} \frac{1}{n^2} \bigg) \notag \\
&= \sup_{n\in\N} \frac{\e^n}{4 n^2} = \infty.
\end{align}
So, by the admissibility theorem above (Theorem~\ref{thm:char infinite-time adm}), $B$ is not infinite-time admissible for $A'$, as desired. $\blacktriangleleft$
\end{ex}

\section*{Acknowledgements}

I would like to thank the German Research Foundation (DFG) financial support through the grant ``Input-to-state stability and stabilization of distributed-parameter systems'' (DA 767/7-1).

\begin{small}

\end{small}

\end{document}